\documentclass[paper,twoside]{article}%%\documentclass[a4paper,twoside]{article}
\usepackage{amsfonts,amsmath,amssymb}

\newtheorem{Definition}{Definition}
\newtheorem{Proposition}{Proposition}
\newtheorem{Theorem}{Theorem}
\newtheorem{Corollary}{Corollary}
\newtheorem{Lemma}{Lemma}

\newenvironment{Proof}{\hfill
\par
\noindent{\bf Proof:}}{\hfill\rule{2mm}{2mm}
\par
\noindent}
\newenvironment{Remark}{\vspace{0.1cm}\noindent{\bf Remark.}}{\vspace{0.1cm}}

\def\vphi{\varphi}
\def\eps{\varepsilon}

\def\d{\mathrm{d}}

\def\R{\mathbb{R}}

\def\AA{\mathcal{A}}

\def\LL{\mathcal{L}}
\def\BB{\mathcal{B}}

\def\DD{\mathcal{D}}
\def\UU{\mathcal{U}}
\def\HH{\mathcal{H}}
\def\VV{\mathcal{V}}

\def\m{\mathrm{m}}

\def\bsob{{\mathrm{H}}_0}
\def\sob{{\mathrm{H}}}

\def\diver{\mathrm{div}}
\def\grad{\mathrm{grad}}

\def\spec{\mathrm{spec}}
\def\bd{\vert_{\partial L(1)}}

\begin{document}

\flushbottom

\begin{center} {\bf \Large Smooth Homogenization of Heat Equations on Tubular Neighborhoods}
\end{center}

\vskip0.5cm

\begin{center} {
O. WITTICH \\\small\em Department of Mathematics and Computer Science, \\ Eindhoven University of Technology, P.O. Box 513, \\ 5600 MB Eindhoven, The Netherlands}
\end{center}

\vskip0.5cm

\begin{quote}
{\bf Abstract.}{\sl We consider the heat equation with Dirichlet boundary conditions on the tubular neighborhood of a closed Riemannian submanifold of a Riemannian manifold. We show that, as the tube diameter tends to zero, a suitably rescaled and renormalized semigroup converges to a limit semigroup in Sobolev spaces of arbitrarily large Sobolev index.}
\end{quote}

\vskip0.5cm

\section{Introduction}

This is the second paper where we investigate the asymptotic behavior of solutions to the heat equation on tubular neighborhoods where Dirichlet boundary conditions are imposed at the boundary of the tubes. As in \cite{pap1}, we consider an $l$-dimensional closed Riemannian manifold $L$ which is isometrically embedded into another $m$-dimensional Riemannian manifold $M$, $m > l$. $L(\eps)$ denotes the tubular neighborhood of radius $\eps > 0$ and the solutions of the heat equation to square - integrable initial conditions are given by the semigroup generated by the Dirichlet Laplacian on $L(\eps)$. If these semigroups are properly rescaled and renormalized, we observe convergence to a limit semigroup as $\eps$ tends to zero. In \cite{pap1}, we proved the statement Theorem \ref{MainSemi1} below about convergence of the semigroups in some suitable $L^2$-space and identified the limit.

In this paper, we will prove Proposition \ref{compactness} below, from which we can immediately conclude Theorem \ref{MainSemi2} which states that the semigroups even converge in Sobolev spaces of arbitrarily large Sobolev index. Consequently, the solutions of the heat equations are smooth functions which converge uniformly with all their derivatives as $\eps$ tends to zero. In the sequel, this behaviour is called {\em smooth convergence}. Note that this is only true for compact subsets of the time axis which do not  contain $t=0$ because neither the semigroups nor the limit $\eps$ to zero are in general continuous at $t=0$.

The result about the smooth convergence of the heat semigroups will be used in \cite{pap3} to complement and to extend results about Brownian motion conditioned to tubular neighborhoods obtained in \cite{SidSmoWeiWit:04} and \cite{SmoWeiWit:07}.

In Section \ref{DLT}, we recall the setup, the construction of the rescaled and renormalized generator family and the main result from \cite{pap1}. We formulate the compactness result Proposition \ref{compactness} and show how we obtain from it Theorem \ref{MainSemi2} about smooth convergence of the semigroups.

In Section \ref{RefLap}, we investigate spectral properties of the Laplacian associated to the reference metric by using results from \cite{BerBou:82} for Riemannian submersions with totally geodesic fibers. This is due to the fact that the projection $\pi : L(1)\to L$ is such a submersion if $L(1)$ is equipped with the reference metric. Finally, we recall the Shapiro-Lopatinskij condition for regular elliptic boundary problems and show how we can equivalently represent the Sobolev norms by the operator norm of such a boundary value problem. This may also serve as a simpler example for what we will do in Section \ref{smooth}.

We think of the generators of the induced dynamics as a perturbation of the generators associated to $g_0$. In this last section, we therefore start by proving some Kato type inequalities relating the perturbation, i.e. the difference of the Laplace operators associated to the induced and to the reference metric, to the reference family $\Delta_0(\eps)$. Then we take up the idea of Section \ref{RefLap} and represent in Proposition \ref{LSindu} the Sobolev norm on the domain of powers of the Dirichlet Laplacian associated to the induced metric by the graph norm of powers of the Laplacian associated to the reference metric. To do so, we have to perform a detailed analysis of the $\eps$-dependence of the respective boundary conditions in Proposition \ref{BigRand} and to use elliptic regularity. By interpolation, we finally arrive at an estimate of the Sobolev norm from above by the graph norm of the induced family $\Delta (\eps)$. This is sufficient to finally prove Proposition \ref{compactness} by an application of the spectral theorem.

\section{The Dirichlet Laplacian on small tubes}\label{DLT}

We consider the {\em Dirichlet Laplacian} $\Delta_{\eps}$ on the Hilbert space $L^2(L(\eps),g)$, i.e. the operator associated to the quadratic form
\begin{equation}
\label{unskaliert}
    q_{(\eps)} (u) := \int_{L(\eps)} \m (dp) \,\Vert \d u\Vert^2
\end{equation}
with domain $\DD (q_{\eps}) = \bsob^1 (L(\eps),g)$. Here, $\m$ denotes the {\em Riemannian volume measure} associated to
$g$ and $\Vert - \Vert$ the norm induced by $g$.

Let $NL\subset TM\vert_L$ be the {\em normal bundle} of the submanifold $L\subset M$. The {\em exponential map} $\exp^M : TM \to M$ restricted to $NL$
thus yields a smooth map $\exp^{\perp} : NL
\to M$. By compactness of $L$, there is some $r > 0$, the {\em injectivity radius}, such that $\exp^{\perp} : U_{\eps}(0)\to
L(\eps)$ is a diffeomorphism from an open $\eps$-neighborhood of the {\em zero section} in $NL$ to the {\em tubular
$\eps$-neighborhood} $L(\eps)$ of $L$ for all $\eps < r$. In the sequel, we assume for simplicity that $r>1$. The {\em
bundle projection} $\pi_N : NL \to L$ induces a submersion $\pi =
\vphi\circ \pi_N \circ\exp^{\perp\,-1} = \exp^{\perp}\circ \pi_N \circ\exp^{\perp\,-1} : L(1) \to L$. For $q\in L$, the preimages $\pi^{-1}
(q)$ provide a decomposition of $L(1)$ into relatively closed, $m-l$-dimensional submanifolds.

On the tube, we consider the family of diffeomorphisms $\sigma_{\eps}:L(\eps)\to L(1)$ given by
\begin{equation}
\label{Reskalierung}
    \sigma_{\eps} (p) := \exp^{\perp} ( \eps^{-1} \,\exp^{\perp\,-1} (p)),
\end{equation}
$\eps > 0$. These maps are called {\em rescaling maps}.

Apart from the metric $g$ induced from $M$, we will also consider another metric, the reference metric $g_0$ on the tubular neighborhood $L(\eps)$. This is the metric induced by the {\em Sasaki metric} (cf. \cite{Sak:97}, (4.6), p. 55 and (4.11), p. 58) on $NL$ via the exponential map.

\begin{Definition}\label{Referenzmetrik} Let $g_S$ be the {\em Sasaki metric} on the total space $NL$ of the normal bundle making it a Riemannian manifold. Via $\exp^{\perp}$, the restriction of the Sasaki metric to $U_1(0)$ is transported to a Riemannian metric on $L(1)$ which we denote by $g_0$. In the sequel, $g_0$ will be referred to as the
{\em reference metric}.
\end{Definition}

Denote the Riemannian volume measure associated to the reference metric by $\m_0$. Denote by
\begin{equation}
\label{Dichte}
    \rho := \frac{\d\m}{\d\m_0} \in C^{\infty}(\overline{L(1)})
\end{equation}
the {\em Radon-Nikodym density} of the volume measures associated to $g$, $g_0$, respectively. The density $\rho > 0$ is
strictly positive, hence pointwise multiplication of functions $f\in L^2 (L(\eps),g_0)$ by $\rho^{-1/2}$ yields a
unitary isomorphism between $L^2 (L(\eps),g_0)$ and $L^2 (L(\eps),g)$ for all $1\geq
\eps > 0$. This, together with the {\em rescaling map} from (\ref{Reskalierung}) yields a unitary isomorphism
$\Sigma_{\eps} : L^2(L(1),g_0)\rightarrow L^2(L(\eps),g)$ given by
\begin{equation}
\label{UniIso}
  \Sigma_{\eps}  u  := \frac{\sigma_{\eps}^*u}{\sqrt{\rho\eps^{m-l}}}
\end{equation}
where we use the shorthand $\sigma_{\eps}^*u := u\circ\sigma_{\eps}$. By smoothness of $\sqrt{\rho}\in
C^{\infty}(\overline{L(1)})$, the maps $\Sigma_{\eps}$ restrict to isomorphisms of the respective domains.

To finally construct the perturbation problem that we will actually investigate, we will {\em rescale} now the quadratic
forms (\ref{unskaliert}) to quadratic forms defined on $\bsob^1(L(1),g_0)$ and {\em renormalize} them with the help of
the lowest eigenvalue $\lambda_0$ of the Dirichlet laplacian for the flat euclidean ball. The rescaled and renormalized
quadratic form
\begin{equation}\label{functional}
F_{\eps} (u) := q_{(\eps)}(\Sigma_{\eps}u) - \lambda_0/\eps^2 \langle u,u\rangle_{0}
\end{equation}
where $\langle -,-\rangle_{0}$ denotes the scalar product on the Hilbert space $L^2(L(1),g_{0})$ associated to the
reference metric, is the quadratic form associated to the operator $$
\Delta (\eps) := \Sigma_{\eps}^{-1}\Delta_{\eps}\Sigma_{\eps} - \lambda_0/\eps^2 . $$

\begin{Definition}\label{deltaeps} For $\eps > 0$, the {\em rescaled family} is given by
\begin{equation*}
\Delta (\eps):= \Sigma_{\eps}^{-1}\,(\Delta_{\eps}-\lambda_{0}/\eps^2)\,\Sigma_{\eps}
\end{equation*} where $\Sigma_{\eps}$ is the unitary map from (\ref{UniIso}) and $\lambda_{0} > 0$ is the lowest
eigenvalue for the Dirichlet boundary problem on the euclidean unit ball $B\subset \R^{m-l}$. the domain of the rescaled family is $\DD (\Delta (\eps))= \bsob^1\cap\sob^2 (L(1),g_{0})$.
\end{Definition}

Consider now again the Dirichlet problem on the flat euclidean ball $B\subset \R^{m-l}$. The eigenspace to the lowest
eigenvalue $\lambda_0 > 0$ is one-dimensional and the eigenfunctions are orthogonally invariant. Let thus $U_0(\vert w
\vert)$ denote a normalized eigenfunction of the flat Dirichlet Laplacian with eigenvalue $\lambda_0$. On the tube, we
consider the function $u_0\in C^{\infty}(\overline{L(1)})$ given by
\begin{equation}
\label{vakuum} u_0(p) := U_0(d(p,\pi(p))) .
\end{equation}

\begin{Definition}\label{EE} Let $g_L$ denote the induced metric on the submanifold $L$. Then $E_{0}\subset L^2
(L(1),g_0)$ is given by
\begin{equation*} E_{0} := \lbrace u \in L^2(L(1),g_0) \,:\, u(p) = u_0(p)\, v(\pi(p)), v\in L^2(L,g_L)\rbrace .
\end{equation*}
\end{Definition}

{\bf Notation.} {\rm (i)} In the sequel, we will denote $E_{0}$ and the $L^2(L(1),g_{0})$-orthogonal projection onto it
by the same symbol. {\rm (ii)} For a given $f\in L^2(L,g_{L})$, we will use the shorthand $\overline{f}=f\circ\pi$ and use the common notion {\em basic functions} for $\overline{f}$ in the sequel.

The following main result of \cite{pap1} is the starting point of our considerations.

\begin{Theorem}\label{MainSemi1} Let $u_{\eps}$, $\eps > 0$ be a sequence in the Hilbert space $L^2(L(1),g_0)$ which converges to $u\in L^2(L(1),g_0)$. Then
$$
\lim_{\eps\to 0}e^{-\frac{t}{2}\Delta (\eps)}u_{\eps}  = E_0\, e^{-\frac{t}{2}(\Delta_L + W_L)} \, E_0 u
$$
strongly in $L^2 (L(1),g_0)$ and uniformly for $t\in K$ where $K\subset (0,\infty)$ is compact. here, $\Delta_L$ denotes the Laplace - Beltrami operator on $L$ and $W_L\in C^{\infty}(L)$.
\end{Theorem}

\begin{Proof} see \cite{pap1}, Theorem 2.
\end{Proof}

By \cite{pap1}, Proposition 1, the {\em effective potential} $W_{L}$ is the restriction of the potential
\begin{equation}\label{effective}
W = \frac{1}{2}\Delta \log\rho -\frac{1}{4}\Vert \d\log\rho\Vert^2 \in C^{\infty}(\overline{L(1)})
\end{equation}
to the submanifold $L$. For an expression of the effective potential in terms of in- and extrinsic geometric quantities such as the scalar curvature of the submanifold and the tension vector field of the embedding,
we refer to \cite{pap1}, Corollary 4. These considerations will play no role in the sequel.

We now want to prove that this convergence result also holds true in a considerably stronger topology. Namely, the
sequence converges in Sobolev spaces of arbitrary large Sobolev index. To show this, we combine the previous result
Theorem \ref{MainSemi1} with the following statement on uniform boundedness which by Rellich's Lemma is a compactness
result on Sobolev spaces, too.

\begin{Proposition}\label{compactness}{\bf(Compactness)} Let $\Vert - \Vert_{k}$ denote the norm of the Sobolev space $\sob^k(L(1),g_0)$. Let $u_{\eps}$, $\eps > 0$ be a sequence of functions in $L^2(L(1),g_0)$ converging in $L^2(L(1),g_0)$ to some $u$. Let $K\subset (0,\infty)$ be compact and
$$
u(\eps,t) := e^{-\frac{t}{2}\Delta (\eps)}u_{\eps}.$$
Then, for all $n\geq 0$ there are numbers $C_n > 0$, $\eps_n > 0$
such that
\begin{enumerate}
\item The sets $\AA_t := \lbrace u(\eps,t) \,:\,\eps \leq \eps_n\rbrace \subset H^{2n}(L(1),g_0)$ are uniformly norm-bounded for all $t\in K$, i.e.
$\Vert u(\eps,t) \Vert_{2n} \leq C_n < \infty$
uniformly for all $\eps \leq \eps_n$.
\item The set $$\AA := \lbrace u(\eps,-)\in C(K,H^{2n}(L(1),g_0)\,:\, \eps \leq \eps_n\rbrace \subset C(K,H^{2n}(L(1),g_0)$$ is uniformly equicontinuous.
\end{enumerate}
\end{Proposition}

From this proposition, we can immediately conclude the main result of this paper.

\begin{Theorem}\label{MainSemi2} For every sequence $u_{\eps}$, $\eps > 0$ with $\lim_{\eps\to 0}u_{\eps}= u$ in $L^2(L(1),g_0)$, we have $$
\lim_{\eps\to 0}e^{-\frac{t}{2}\Delta (\eps)}u_{\eps}  = E_0\, e^{-\frac{t}{2}(\Delta_L + W_L)} \, E_0 u $$
strongly in every Sobolev space $\sob^{k} (L(1),g_0)$, $k\geq 0$, and uniformly for $t\in K$ where $K\subset (0,\infty)$ is compact.
\end{Theorem}

\begin{Proof} By Proposition \ref{compactness} {\rm(i)} and by compactness of the embeddings $$\sob^{2k}(L(1),g_{0})\subset
\sob^{2k-1}(L(1),g_{0},$$ the sets $\AA_t\subset
\sob^{2k-1}(L(1),g_{0})$ are relatively compact. By Proposition \ref{compactness}, {\rm (ii)} and the Arzela - Ascoli theorem that implies that
$\AA \subset\subset C(K,\sob^{2k-1}(L(1),g_{0}))$ is relatively compact. We consider
now an arbitrary sequence $a_{n,n\geq 1}$ in $\AA$. By compactness, it contains a subsequence $a'_n,n\geq 1$ which
converges to an element  $a'_{{\infty}}\in C(K,\sob^{2k-1}(L(1),g_{0}))$. Then $$
\lim_{n\to\infty} \sup_{t\in K}\Vert a'_{n}(t) - a'_{\infty}(t) \Vert_{0} \leq \lim_{n\to\infty} \sup_{t\in K}\Vert a'_{n}(t) - a'_{\infty}(t) \Vert_{2k-1}
= 0 $$
and hence, the subsequence also converges in $C(K,L^2(L(1),g_{0}))$ to $a'_{\infty}$. By Theorem \ref{MainSemi1}, that
implies that $$
a'_{\infty}(t)= E_{0}\,e^{-\frac{t}{2}(\Delta_{L}+W_L)}\,E_{0}u.$$
Hence, every subsequence contains a convergent subsequence and all these convergent subsequences have the same limit.
That implies the statement.
\end{Proof}

In the remainder, all we have to do is to prove Proposition \ref{compactness}.

\section{Reference metric and associated Laplacian}\label{RefLap}

In this section, we investigate some properties of the reference metric and the associated Laplacian. The basic idea of the proof of the homogenization result is to compare the situation for the induced metric with
the situation for the reference metric. Thus, it is natural to introduce as well a rescaled family for the Laplacian
associated to $g_{0}$. Thus, let $\Delta_{0}:=-\diver_{g_{0}}\grad$ denote the Laplace - Beltrami operator on $L(1)$
with respect to the metric $g_{0}$. The corresponding {\em Dirichlet Laplacian} $\Delta_{0,\eps}$ is given by
$\Delta_{0,\eps}u := \Delta_{0} u$ for all $u$ in the domain $\DD (\Delta_{0,\eps}):=\bsob^1\cap\sob^2 (L(\eps),g_{0})$.
Thus, $\Delta_{0,\eps}$ is the positive operator associated to the quadratic form $$
q_{0,\eps}(u):=\int_{L(\eps)} \m_{0}(dp) \Vert \d u
\Vert^2_{0} $$
with domain $\DD (q_{0,\eps})=\bsob^1 (L(\eps),g_{0})$. Specializing the construction from Section \ref{DLT}, (\ref{UniIso}), to the
reference metric basically means that due to $\rho \equiv 1$, we obtain  $\Sigma_{0,\eps}u =
\eps^{(l-m)/2}\sigma_{\eps}^{*}u$.

\begin{Definition}\label{deltanulleps} For $\eps > 0$, the {\em rescaled reference family} is given by
\begin{equation*}
\Delta_{0}(\eps):= \Sigma_{0,\eps}^{-1}\,(\Delta_{0,\eps}-\lambda_{0}/\eps^2)\,\Sigma_{0,\eps} =
\sigma_{\eps}^{-1}\,(\Delta_{0,\eps}-\lambda_{0}/\eps^2)\,\sigma_{\eps}
\end{equation*} where $\sigma_{\eps}$ is the rescaling map from (\ref{Reskalierung}) and $\lambda_{0} > 0$ is the lowest
eigenvalue for the Dirichlet boundary problem on the euclidean unit ball $B\subset \R^{m-l}$. the domain of the rescaled reference family is $\DD (\Delta_{0}(\eps))= \bsob^1\cap\sob^2 (L(1),g_{0})$.
\end{Definition}

The rescaled reference family can be quite well understood due to special geometric properties of the reference metric.
It turns out that with respect to the reference metric, the canonical projection is a {\em riemannian submersion} with
{\em totally geodesic fibers}. That implies that the Laplacian on $L(1)$ can be decomposed into a {\em vertical} and a
{\em horizontal} part and that these two parts actually commute. This property provides us with quite detailed
information about the spectrum of the Laplacian and its behaviour under rescaling which makes the reference metric
especially suitable to serve as a basis for our perturbational ansatz. In a situation where the fibers are closed and,
consequently, no boundary conditions are present, these properties are proved in \cite{BerBou:82}. In the sequel, we
will basically follow this exposition with minor modifications due to the presence of Dirichlet boundary conditions.

\subsection{The decomposition of the Laplacian}\label{BerBou}

Let $\phi:M\rightarrow N$ denote a map between Riemannian manifolds and denote by $\phi_{*}:TM\rightarrow TN$ its {\em
tangent map}. For $p\in  M$, the space $\VV_{p}\subset T_{p}M$ of {\em vertical vectors} is given by $\VV_{p} :=\ker
\phi_{*\,p}$. The space $\HH_{p}\subset T_{p}L(1)$ of {\em horizontal vectors} is given by the orthogonal complement
$\HH_{p}(1):=\ker\phi_{*\,p}^{\perp}$ in $T_{p}M$. In the same way, we can define horizontal and vertical subbundle
$\HH,\VV\subset TM$ and horizontal and vertical vector fields. Clearly, $\HH\oplus \VV = TM$. Assuming that the fibers
$\phi^{-1}(q)\subset M$ are Riemannian submanifolds for all $q\in N$, we may consider the Laplace-Beltrami operators
$D_{q}$ on the fibers $\phi^{-1}(q)$ with respect to the induced metric.

\begin{Definition}\label{horivert} Let $f\in C^2 (M)$. {\rm (i)} The operator $$
\Delta_{M}^Vf (p):= D_{\phi (p)}\left(f\vert_{\phi^{-1}(\phi(p))}\right)(p)$$
is called the {\em vertical Laplacian} on $M$. {\rm (ii)} The operator $$
\Delta_{M}^Hf:= (\Delta_{M}-\Delta_{M}^V)f$$
is called the {\em horizontal Laplacian} on $M$.
\end{Definition}

$\phi$ is called a {\em Riemannian submersion}, iff it leaves the norm of horizontal vectors unaffected, i.e.
$\Vert\phi_{*}v\Vert_{p} = \Vert v \Vert_{\phi(p)}$ for all $v\in \HH_{p}$, $p\in M$. The map is said to have {\em totally geodesic fibers}, iff the fibers
$\phi^{-1}(q)\subset M$ are totally geodesic submanifolds (see \cite{KobNom:63}, p. 180).\\

\begin{Theorem}\label{comm1} Let $f\in C^4(M)$. If $\phi : M\to N$ is a Riemannian submersion with totally
geodesic fibers, we have $
\lbrack \Delta_{M},\Delta_{M}^V \rbrack f = \lbrack \Delta_{M},\Delta_{M}^H \rbrack f = \lbrack
\Delta_{M}^H,\Delta_{M}^V \rbrack f = 0$.
\end{Theorem}

\begin{Proof} \cite{BerBou:82}, Theorem 1.5.
\end{Proof}

We consider the projection map $\pi : L(1)\rightarrow L$. Actually, we are in the situation considered
above. In particular, the operators $\Delta_{0}$, $\Delta_{0}^V$ and $\Delta_{0}^H$ commute in the sense just described.

\begin{Lemma}\label{PiRiemSub} If $L(1)$ is equipped with the reference metric $g_{0}$, the canonical projection $\pi :
L(1)\rightarrow L$ is a Riemannian submersion with totally geodesic fibers.
\end{Lemma}

\begin{Proof} By \cite{pap1}, Section 6.1, we can find local coordinates $(x,w)$ where $x = (x^i)_{i=1,...,L}$ and $w = (w^{\alpha})_{\alpha = 1,...,m-l}$ for the tube such that
$\pi (x,w) = x$, the horizontal spaces $\HH_{x,w}$ are generated by the vectors
\begin{equation}\label{horizont}
X_i (x,w) = \partial_i - w^{\alpha}C_{i\alpha}^{\beta}\partial_{\beta}
\end{equation}
and the reference metric is given by
\begin{equation}\label{Referenzmetrik lokal}
g_0(x,w)  =  \left(
\begin{array}{c|c} g_{L,ij} + w^{\mu}w^{\nu}C_{i\mu}^{\alpha}C_{j\nu}^{\beta}\delta_{\alpha\beta} & w^{\mu}C_{i\mu}^{\alpha} \\
\hline w^{\nu}C_{j\nu}^{\beta} & \delta_{\alpha\beta}
\end{array}\right)\end{equation}
where $C_{i\alpha}^{\mu}$ denote the connection coefficients of the bundle connection on $NL$. By
$$
\langle X_i, X_j\rangle_0 = g_{L,ij} = \langle \partial_i,\partial_j\rangle_L =  \langle \pi_* X_i,\pi_*X_j\rangle_L ,
$$
we have a Riemannian submersion. The vanishing of the second fundamental form of the fibers follows from the vanishing of the Christoffel symbols $\Gamma^i_{\alpha\beta}=0$
and $\Gamma^{\mu}_{\alpha\beta}=0$ for $g_0$.
\end{Proof}

Thus, we finally obtain by combining the two statements above

\begin{Corollary}\label{commutative} Let $f\in C^4(\overline{L(1)})$. Then we have $$
\lbrack \Delta_{0},\Delta_{0}^V \rbrack f = \lbrack \Delta_{0},\Delta_{0}^H \rbrack f = \lbrack
\Delta_{0}^H,\Delta_{0}^V \rbrack f = 0.$$
\end{Corollary}

\subsection{The vertical operator as direct integral}\label{VODI}

In order to make use of the preceding discussion also for the investigation of the spectral properties of the Dirichlet
Laplacian as an operator on $L^2(L(1),g_{0})$, we consider the tube $L(1)$ equipped with the reference metric $g_0$.
First, we will construct a direct integral decomposition of $L^2 (L(1),g_0)$ along the leaves and a decomposable
self-adjoint operator which extends the vertical Laplacian defined for functions $f\in C(\overline{L(1)})\cap C^2(L(1))$
with Dirichlet boundary conditions $f\vert_{\partial L(1)} = 0$.

Now we proceed by constructing the leaf operator as a direct integral of Dirichlet operators on the fibers. Let
$D_{q,q\in L}$ be given by the family of {\em Dirichlet operators} on the leaves. To be precise let $D_q := -\diver_q
\grad_q$, the Laplace-Beltrami operator on $(\pi^{-1}(q))$ with domain
\begin{equation*}
\DD (D_q) := \bsob^1\cap\sob^{2} (\pi^{-1}(q)).
\end{equation*} Since the differential expression for the Laplacian depends smoothly on the fibers, the family is
measurable. Since for riemannian submersions with totally geodesic fibers the fibers are isometric (isometries are
generated by horizontal vector fields, see
\cite{BerBou:82}), we check by a calculation in Fermi coordinates that all fibers can
be isometrically mapped to the Euclidean unit ball $B\subset
\R^{m-l}$. Hence, the operators $D_q$ are unitarily equivalent to the Dirichlet Laplacian on $B\subset \R^{m-l}$ for all
$q\in L$. Hence, they are self-adjoint, and their spectra $\spec (D_q)$ are {\em semi-simple} (i.e. consist only of
eigenvalues of finite multiplicity) and strictly positive with smallest eigenvalue $\lambda_0 > 0$. Hence, the operator
\begin{equation}
\Delta_{0}^V := \int_L^{\oplus} \m_L (dq) \, D_q
\end{equation}
is self-adjoint and extends the operator $\Delta_{0}^V$ defined on functions $f\in C(\overline{L(1)})\cap C^2(L(1))$
with Dirichlet boundary conditions $f\vert_{
\partial L(1)} = 0$. It is therefore justified to denote it by the same symbol. Hence, the Dirichlet leaf operator
\begin{equation}
\label{DLO} \Delta_{0}^V - \lambda_0 := \int_L^{\oplus} \m_L (dq) \, (D_q -
\lambda_0)
\end{equation}
is as well self-adjoint and $$
\spec (\Delta_{0}^V-\lambda_{0}) := \lbrace \lambda_{k}-\lambda_{0}\,:\, k\geq 1\rbrace$$
where $\lambda_{k,k\geq 0}$ denotes the spectrum of the Dirichlet Laplacian on the Euclidean unit ball
$B\subset\R^{m-l}$. As already said in the introduction, the eigenspace of the Dirichlet Laplacian in $B\subset\R^{m-l}$
belonging to the smallest eigenvalue $\lambda_0$ is one-dimensional, and the eigenfunctions are orthogonally invariant.
Let again $u_0 (w) := U_0 (\vert w
\vert)$ be a normalized eigenfunction generating this eigenspace. Recall, \cite{ReeSim:78}, XIII.16, Definition, p. 281,
that a bounded operator is called {\em decomposable}, if it can be written as a direct integral of operators as above.
Then the kernel coincides with the asymptotic subspace from Definition
\ref{EE}.

\begin{Proposition}\label{kern Dirichlet} {\rm (i)} The Dirichlet leaf operator $\Delta_{0}^V - \lambda_0$ is
non-negative with kernel
\begin{equation*}
\ker (\Delta_{0}^V - \lambda_0) = \lbrace u\in V\, :\, u(p):= u_0 (p)\, v(\pi (p)), v\in L^2 (L,g_L)\rbrace.
\end{equation*} Here, $u_0 (p) := U_0(d_0(p,\pi(p)))$, where $d_0$ denotes the Riemannian distance with respect to the
reference metric $g_0$ on $L(1)$. {\rm (ii)} The orthogonal projection $E_{0}$ onto $\ker (\Delta_{0}^V - \lambda_0)$ is
decomposable and given by T
\begin{equation*} E_{0} = \int_L^{\oplus} \m_L (dq) E_{0,q} .
\end{equation*}
\end{Proposition}

\begin{Remark} {\rm (i)} Note that in local Fermi coordinates $p\equiv (x,w)$ we have $$
u_0(p) = U_{0}(d(p,\pi(p))) = U_0(\vert w\vert) = u_0(w) $$
which justifies our slight abuse of notation. {\rm (ii)} In particular, the corollary implies that the kernel coincides
with the asymptotic subspace $E_0$ from Definition \ref{EE}.
\end{Remark}

\begin{Proof} see \cite{pap1}, Corollary 5.
\end{Proof}

A fortiori, by \cite{ReeSim:78}, Theorem XIII.85, p. 284, we may conclude analogous results for all eigenspaces. We thus
obtain the following spectral decomposition of the operator $\Delta_{0}^V$: Let $\lambda_{k,k\geq 0} $ be the collection
of eigenvalues of the Dirichlet Laplacian on the euclidean unit ball $B\subset \R^{m-l}$ as above and $P_{k,k\geq 1}$
together with $P_{0}=E_{0}$ denote the corresponding eigenspaces. Again, we denote the eigenspaces and the corresponding
orthogonal projections by the same symbol. Then, by mapping the euclidean unit ball $B$ isometrically onto the fiber
$\pi^{-1}(q)$, we obtain projections $E_{k,q}$ on $L^2(\pi^{-1}(q),\m_{q})$ induced by $P_{k}$. Then,  by
\cite{ReeSim:78}, Theorem XIII.85, p. 284, we can compute the spectral decomposition of $D$.

\begin{Proposition}\label{Dspektral} The operator $\Delta_{0}^V$ is self - adjoint on $L^2(L(1),g_0)$ with spectral
decomposition
\begin{equation*}\label{Dspektrum} \Delta_{0}^V = \lambda_{0} E_{0} + \sum_{k\geq 1}\lambda_{k} E_{k}
\end{equation*} where
\begin{equation*} E_{k} = \int_L^{\oplus} \m_L (dq) E_{k,q}.
\end{equation*}
\end{Proposition}

\begin{Proof} see \cite{pap1}, Proposition 3.
\end{Proof}

\begin{Remark} Note that the spectral eigenspaces $E_{k}$ are infinite-dimensional in general. In the next subsection,
we will show that they all consist of eigenfunctions of the Dirichlet Laplacian $\Delta_{0}$.
\end{Remark}

\subsection{Common eigenfunctions and spectra}\label{Cefu}

In this subsection, we show that the self-adjoint Dirichlet Laplacian $\Delta_{0}$ with domain
$\bsob^1\cap\sob^2 (L(1),g_{0})$ actually commutes with the direct integral operator just defined, which is a
self-adjoint extension of $\Delta_{0}^V$ and which we already decided to denote by the same symbol. We start from the
following well known result (for a proof see for instance \cite{Tay:99}, Ch. 5.1, p. 303 ff.).

\begin{Theorem}\label{SpekDirichlet} Let $M$ be a smooth riemannian manifold with smooth boundary $\partial M$. The
Dirichlet Laplacian is self-adjoint on $L^2(M,g)$ with a spectrum only consisting of eigenvalues of finite multiplicity.
The corresponding eigenfunctions are smooth and can be chosen to form an orthonormal base of $L^2(M,g)$.
\end{Theorem}

\begin{Lemma} {\rm (i)} Let $u$ be an eigenfunction of $\Delta_{0}$ with corresponding eigenvalue $\mu$. Then
$\Delta_{0}^Vu$ is as well an eigenfunction of $\Delta_{0}$ with the same eigenvalue. {\rm (ii)} $\Delta_{0}^V:E\to E$
yields a self-adjoint isomorphism for every eigenspace $E$ of $\Delta_{0}$.
\end{Lemma}

\begin{Proof} {\rm (i)} By Theorem \ref{SpekDirichlet}, $u$ is smooth and we may apply Corollary \ref{commutative}. Hence
$$
\Delta_{0}(\Delta_{0}^Vu) = \Delta_{0}^V(\Delta_{0}u) = \mu\Delta_{0}^Vu. $$
{\rm (ii)} By Theorem \ref{SpekDirichlet}, the eigenspaces are finite-dimensional. By $\Delta_{0}^V > 0$, the mapping is
injective and hence surjective, too. Thus, $\Delta_{0}^V$ is self-adjoint on $L^2(L(1),g_{0})$ and leaves $E$ invariant.
Hence it is also self-adjoint when restricted to $E$.
\end{Proof}

By the {\em spectral theorem} we now diagonalize $\Delta_{0}^V$ on the respective eigenspaces of $\Delta_{0}$. Thus, we
find smooth common eigenfunctions for the Laplacian and the vertical Laplacian. Since these functions are already known
to form an orthonormal base of $L^2(L(1),g_{0})$ by Theorem \ref{SpekDirichlet}, they also provide a spectral
decomposition of the self-adjoint direct integral operator representing the vertical Laplacian. Thus, we obtain the
following almost complete spectral characterization of the three operators involved.

\begin{Proposition}\label{commonEF} Let the eigenvalues $\mu_{0}\leq\mu_{1}\leq...$ of the Dirichlet Laplacian
$\Delta_{0}$ be ordered by magnitude. Then there is an orthonormal basis $u_{s,s\geq 1}$ of $L^2(L(1),g_{0})$ consisting
only of smooth functions such that
\begin{equation*}
\begin{array}{lll}
{\mathrm{(i)}}\,\, \Delta_{0}u_{s}= \mu_{s}u_{s},&
{\mathrm{(ii)}}\,\,\Delta_{0}^Vu_{s}=\lambda_{k(s)}u_{s},&
{\mathrm{(iii)}}\,\,\Delta_{0}^Hu_{s}=(\mu_{s} - \lambda_{k(s)})u_{s},
\end{array}
\end{equation*}
where $k(s)$ is given by the unique $k\in\lbrace 0,1,...\rbrace$ such that $E_{k(s)}u_{s}=u_{s}$.
\end{Proposition}

The spectrum of $\Delta_{0}^H$ hence consists of the differences of eigenvalues of $\Delta_{0}$ and $\Delta_{0}^V$.
However, different eigenvalues may lead to the same differences which makes it difficult to make definite conclusions
about the spectrum of $\Delta_{0}^H$. It might even be that the spectrum is no longer discrete (cf. \cite{BerBou:82},
3.4). However, note that by (\ref{Referenzmetrik lokal}), in local coordinates the quadratic form associated to $\Delta_0^H$ is given by
\begin{equation}\label{horipositiv}
\tau (u) = \int dxdw\,\sqrt{\det g_L}\, g_{L}^{ij}\, X_iu\, X_ju
\end{equation}
where $X_i $ are the horizontal vector fields from (\ref{horizont}). For a proof see \cite{pap1}, Lemma 4. In particular, by the positive definiteness of $g_L$, we have $\tau (u) \geq 0$.

\subsection{The rescaled reference family}\label{viervier}

With these preliminaries, we can now derive two properties of the rescaled reference family. This is the starting
point to consider the full rescaled reference family as a singular perturbation of this expression.\\

Under the rescaling map $\sigma_{\eps}$, the reference metric scales as
\begin{equation}\label{can_variation}
g_0(\eps) = \sigma_{\eps}^{*\,-1}g_0 =  \left(
\begin{array}{c|c} g_{L,ij} + w^{\mu}w^{\nu}C_{i\mu}^{\alpha}C_{j\nu}^{\beta}\delta_{\alpha\beta} & w^{\mu}C_{i\mu}^{\alpha} \\
\hline w^{\nu}C_{j\nu}^{\beta} &\eps^2\, \delta_{\alpha\beta}
\end{array}\right)\end{equation}
in the local coordinates considered already in (\ref{Referenzmetrik lokal}). Hence, only the fiber metric is affected by the rescaling. In
\cite{BerBou:82}, Sec. 5, this scaling behaviour is called the {\em canonical variation} of a Riemannian submersion. If
$L(1)$ is equipped with any of the metrics $g_{0}(\eps)$, $\eps > 0$, $\pi$ is still a Riemannian submersion with
totally geodesic fibers (see \cite{BerBou:82}, Prop. 5.2). We start with the following decomposition formula (cf.
\cite{BerBou:82}, Prop. 5.3).

\begin{Proposition}\label{DecoReference} The rescaled family can be written
\begin{equation*}
\Delta_{0}(\eps) = \frac{1}{\eps^2}\left(\Delta_{0}^V-\lambda_{0}\right) + \Delta_{0}^H .
\end{equation*}
\end{Proposition}

\begin{Proof} By \cite{BerBou:82}, Prop. 5.3, the Laplacian on $(L(1),g_0(\eps))$ is given by
$$
\sigma_{\eps}^{-1}\,\Delta_{0,\eps}\,\sigma_{\eps} = \frac{1}{\eps^2} \Delta_{0}^V + \Delta_{0}^H $$
which implies the statement after subtracting $\lambda_{0}/\eps^2$.
\end{Proof}

The operators $\eps^{-2}(\Delta_{0}^V - \lambda_{0})$ and $\Delta_{0}^H$ are still commuting.
The second property is the following lemma about the eigenvalues of the vertical Laplacian which will be used several
times in the sequel.

\begin{Lemma}\label{smoothEV} Let $\lambda_0\leq\lambda_1\leq...$ be the eigenvalues of the Dirichlet Laplacian for the
unit ball $B\subset\R^{m-l}$. Then we have: {\rm (i)} If $\eps \leq 1 - \lambda_0/\lambda_1$, then $(\lambda_k
-\lambda_0)/\eps^2\geq \lambda_k/\eps$ for all $k\geq 1$. {\rm (ii)} If $\eps^2 \leq 1 - \lambda_0/\lambda_1$, then
$(\lambda_k -\lambda_0)/\eps^2\geq \lambda_k$ for all $k\geq 1$.
\end{Lemma}

\begin{Proof} Both parts follow from the inequality $$
\frac{\lambda_k - \lambda_0}{\lambda_{k}}= 1 - \frac{ \lambda_0}{\lambda_{k}}\geq 1 - \frac{
\lambda_0}{\lambda_{1}}=\frac{\lambda_1 - \lambda_0}{\lambda_{1}}$$
and the fact that the lowest eigenspace of the Dirichlet laplacian for the flat unit ball is non-degenerate, i.e $\lambda_1 > \lambda_0$.
\end{Proof}

\subsection{Sobolev norms and the reference Laplacian}\label{Lapnorm}

One cornerstone of the proof is the fact that on several subspaces of the Sobolev spaces $H^{2k}(L(1),g_0)$, the Sobolev
norm can be equivalently expressed by powers of the Laplace-Beltrami operator. Since $L(1)$ is a manifold with smooth
boundary, we will need some properties of associated boundary problems. Let thus $\LL$ be an elliptic operator of order
$2k$ on $L(1)$ with coefficients that are smooth on $\overline{L(1)}$. Let furthermore $\BB := \lbrace \BB_1,...,\BB_k
\rbrace$ be a system smooth differential expressions defined in a neighborhood of $\partial L(1)$ in $M$. The following
definition is adapted from \cite{Agr:97}, Sec. 1.3.

\begin{Definition}{\bf (Shapiro - Lopatinskij condition)} The system $\lbrace \LL, \BB\rbrace$ is said to provide a {\em
regular elliptic boundary problem} iff
\begin{enumerate}
\item $0\leq m_j = \mathrm{ord} \BB_j < 2k$ for $j=1,...,k$ and $m_i\neq m_j$ for $j\neq k$,
\item the {\em principal symbols} $\BB_j^o (x,\xi)$ are such that $\BB_j^o (x,n(x))\neq 0$ for all $x\in\partial L(1)$,
where $n(x)$ is the unit normal vector on the boundary,
\item for any vector $0\neq \xi\in T_x
\partial L(1)$, the polynomials $\BB_j^o(\tau,\xi):=\BB_j^o (x,\xi + \tau n(x))$ of the complex variable $\tau$ are
linearly independent modulo the polynomial $$
\LL_+^o (\tau) := \prod_{j=1}^k (\tau - \tau_j^+(x,\xi)) $$
where the numbers $\tau_j^+(x,\xi)$ are the zeros of the polynomial $\LL^o (x,\xi + \tau n(x))$ with positive imaginary
parts.
\end{enumerate} {\rm (i)} - {\rm (iii)} are called the {\em Shapiro - Lopatinskij condition}.
\end{Definition}

Since the regularity of an elliptic boundary value problem implies the validity of a certain coercivity estimate and,
finally, the subsequent representation of the Sobolev norm, to prove the following statement we basically have to check
the Shapiro-Lopatinskij condition for the boundary problem we have in mind.

Let now $\Delta_0$ denote the Dirichlet Laplacian associated to $g_0$, i.e. with zero boundary conditions on $\partial L(1)$ and denote by $\Delta_0^k$ ist kth power.

\begin{Lemma}\label{normlap} On $\DD (\Delta_0^k)\subset \sob^{2k}(L(1),g_0)$, the $2k$-Sobolev norm is equivalent to
the norm $$
\vert\Vert u \Vert\vert_{2k} := \Vert u \Vert_0 + \Vert \Delta_0^k u \Vert_0 . $$
\end{Lemma}

\begin{Proof} It is sufficient to prove the regularity result
\begin{equation}
\label{estima}
\Vert u \Vert_{2k} \leq C\left(\Vert u \Vert_0 + \Vert \Delta_0^k u\Vert_0 \right)
\end{equation}
since the reversed inequality follows from the fact that $\Delta_0^k : \sob^{2k}(L(1),g_0) \to L^2(L(1),g_0)$ is
continuous. To prove (\ref{estima}), we use that such an estimate always holds for (regular) elliptic boundary problems
(see e.g. the survey \cite{Agr:97}, Theorem 2.2.1, p. 16). Hence, it remains to prove that $\LL = \Delta_0^k$ with
boundary conditions $$
\BB_1u = u\bd = \BB_2 u = \Delta_0u\bd = ... = \BB_k u = \Delta^{k-1}_0 u = 0$$
satisfies the Shapiro - Lopatinskij condition. But for the boundary problem above, we can do all computations
explicitly. We have thus $$
\begin{array}{cc} \LL_+^o (\tau,\xi) = (\tau - i\Vert\xi\Vert)^k, & \BB_l^o (\tau,\xi) = \lbrace(\tau +
i\Vert\xi\Vert)(\tau - i\Vert\xi\Vert)\rbrace^{l-1}
\end{array} $$
where $l=1,...,k$. Let now $a=(a_1,...,a_k)$ and $$
Q_a(\tau,\xi) := \sum_{l=1}^k a_l\, \BB_l^o(\tau,\xi). $$
Suppose that $l_0 := \min \lbrace s \,:\, a_s \neq 0\rbrace$. Then $$
\frac{Q_a(\tau,\xi)}{\LL_+^o (\tau,\xi)} = \frac{Q_a(\tau,\xi)}{(\tau - i\Vert\xi\Vert)^{l_0-1}}\,:\,\frac{\LL_+^o
(\tau,\xi)}{(\tau - i\Vert\xi\Vert)^{l_0-1}} = \frac{Q_{a,l_0}(\tau,\xi)}{\LL_{+,l_0}^o (\tau,\xi)}, $$
$Q_{a,l_0}$ and $\LL_{+,l_0}^o$ are both polynomials, and $\LL_{+,l_0}^o$ has only one zero of multiplicity $k-l+1$ at
$i\Vert\xi\Vert$ whereas $Q_{a,l_0}(i\Vert\xi\Vert,\xi)=(2i\Vert\xi\Vert)^{l_0-1}a_{l_0} \neq 0$. Hence, the
Shapiro-Lopatinskij condition holds.
\end{Proof}

\section{Smooth convergence}\label{smooth}

Now we are going to prove convergence of the semigroups generated by the rescaled Laplacians associated to the induced metric in Sobolev spaces of arbitrarily large Sobolev index (in the sequel called {\em smooth convergence}). By the proof of Theorem \ref{MainSemi2}, we have already reduced that to the proof of Proposition
\ref{compactness} on uniform boundedness of the Sobolev norms. That implies
compactness of the family in Sobolev spaces of arbitrarily large index by Rellich's Lemma (see \cite{Tay:99}, Prop. 4.4,
p. 287). Thus, the family actually converges in these Sobolev spaces and the limit must coincide with the limit in $L^2(L(1),g_0)$ from Theorem \ref{MainSemi1}. The proof of Proposition \ref{compactness} is based on the representation of
the Sobolev norms using the reference Laplacian in Lemma \ref{normlap} and a perturbation ansatz for the boundary problem.

\subsection{The induced metric}

In the sequel, we will make use of the perturbation expansion of the operator family associated to the rescaling of the induced metric in local Fermi
coordinates. Since the Laplacian is a (non linear) functional of the metric, all information about the families is
contained in the perturbation expansion of the difference of induced and reference metric. We summarize the necessary
facts in the following theorem, the three statements of which were proven in Section 5.2 of \cite{pap1}.

\begin{Theorem}\label{Fermi Metrik} {\rm (i)} In local Fermi coordinates the metric tensor is given by
\begin{equation*}
    g(x,w) =
    \left(
\begin{array}{c|c}
    1 & cb^{-1} \\
\hline 0 & 1
    \end{array}\right) \,
    \left(
\begin{array}{c|c}
    a  & 0 \\
\hline 0 & b
    \end{array}\right)\,
    \left(
\begin{array}{c|c}
    1 & 0 \\
\hline b^{-1} c^+ & 1
    \end{array}\right)
\end{equation*} where
\begin{eqnarray*}
a_{ij} (x,w) &=& g_{L,ij} (x) - 2w^{\alpha} g_{L\,is} A_{\alpha j}^s
\\
& & +  w^{\alpha} w^{\beta} ( g_{L\,rs}A_{\alpha i}^r A_{\beta j}^s  -  R_{i\alpha
j\beta}) + O(\vert w
\vert^3),
\\
c_{i\sigma} (x,w) &=& w^{\alpha} C_{i\alpha}^{\sigma} (x) + O(\vert w
\vert^2),
\\
b_{\mu\sigma} (x,w) &=& \delta_{\mu\sigma} -
\frac{1}{3}w^{\alpha}w^{\beta} R_{\mu\alpha\sigma\beta} (x) + O(\vert w \vert^3).
\\
\end{eqnarray*}
Here, $g_L$ denotes the metric on $L$, $A_{\alpha}$
the {\em Weingarten map} of the embedding, $C_{i\alpha}^{\mu} (x) :=
\langle \nu_{\mu}(x), D_i \nu_{\alpha}(x)\rangle$ denote the {\em connection coefficients} of the {\em induced
connection} $D$ on the normal bundle and $R$ is the {\em Riemannian curvature tensor} of $M$.

{\rm (ii)} Let $g(\eps) := \sigma_{\eps}^{-1\,*}g$ be the {\em rescaled metric}. For the
dual metric we have $g^{*}(\eps) = g_{0}^*(\eps) + H^*(\eps)$ where
\begin{equation*}
    \sup_{p\in L(1)} \Vert H^*(\eps) - H^*(0)\Vert = O(\eps)
\end{equation*} as $\eps$ tends to zero. Furthermore, in local coordinates,
\begin{equation*}
    H^*(0)=\frac{1}{3} w^{\alpha} w^{\beta}\left(
\begin{array}{c|c}
    0& 0 \\
\hline
    0 &
     R_{\mu\alpha\nu\beta}
    \end{array}\right).
\end{equation*}

{\rm (iii)} We have $W_{\eps} = W_L\circ\pi + \eps W(\eps)= \overline{W}_L\circ\pi + \eps W(\eps)$ where $W_L$ denotes the restriction of $W$ to $L$ and $W(\eps)$ is a smooth bounded function on $L(1)$ converging uniformly as $\eps$ tends to zero.
\end{Theorem}

\subsection{A decomposition of $\Delta (\eps)$}

By \cite{pap1}, Prop. 1, the rescaled and renormalized quadratic form (\ref{functional}) is given by
$$
F_{\eps}(u) = \int_{L(1)} \m_0(dq)\left(\Vert \d u \Vert_{g(\eps)}^2 + (W_{\eps}- \frac{\lambda_0}{\eps^2})u^2\right)
$$
where $g(\eps) = \sigma_{\eps}^{-1\,*}g$, $W_{\eps} = W \circ \sigma_{\eps}^{-1}$ and $W$ is the effective potential (\ref{effective}). By Theorem \ref{Fermi Metrik}, we can even rewrite this as
$$
F_{\eps}(u) = \int_{L(1)} \m_0(dq)\left(\Vert \d u \Vert_{g_0(\eps)}^2 - \frac{\lambda_0}{\eps^2}u^2 + \langle \d u, \d u\rangle_H + \overline{W}_Lu^2 + \eps r(\eps,u,\d u) \right)
$$
where $r(\eps, u, \d u)$ is a closed quadratic form on $\bsob^1(L(1),g_0)$ with coefficients that converge uniformly as $\eps$ tends to zero. Now the operator family associated to the first term of $F_{\eps}$ is the renormalized reference family. Now, in local coordinates the perturbation is given by
\begin{eqnarray*}
 \langle \d u, \d u\rangle_H &=& \frac{1}{3}\int_{L(1)} \m_0(dq) \,w^{\alpha}w^{\beta}R_{\alpha\mu\beta\nu}\partial_{\mu}u\,\partial_{\nu}u \\
 &=&  \frac{1}{12}\int_{L(1)} \m_0(dq) \,R_{\alpha\mu\beta\nu}\,L_{\mu\alpha}u\,L_{\nu\beta}u \\
 &=&  \frac{1}{12}\int_{L} \m_L(dp) \,R_{\alpha\mu\beta\nu}(p)\,\langle L_{\mu\alpha}u_p, L_{\nu\beta}u_p\rangle_p \end{eqnarray*}
where $L_{\mu\alpha} = w^{\alpha}\partial_{\mu} - w^{\mu}\partial_{\alpha}$ and we used the symmetries of the curvature tensor. That implies that for $u\in\bsob^1\cap\sob^2 (L(1),g_0)$, partial integration leads to
\begin{equation*} F_{\eps}(u)=\frac{1}{2}\left\lbrace  \langle u,\Delta_{0} (\eps)u\rangle_{0} + \langle u, Au
\rangle_{0} + \eps\langle u,R(\eps)u\rangle_{0}\right\rbrace
\end{equation*}
where $A$ is a second order differential operator with smooth coefficients and $R(\eps)$ is a second
order differential operator with uniformly bounded smooth coefficients. Thus, the operator family can be decomposed into
\begin{equation}\label{decomfamily}
\Delta (\eps)u =\Delta_{0}(\eps) u + Au + \eps R(\eps)u
\end{equation}
where
\begin{equation}\label{Alokal}
A (x,w) = -
\frac{1}{12} R_{\beta\nu\mu\alpha}(x)\,L_{\beta\nu}L_{\alpha\mu} + W_{L}(x).
\end{equation}
In the sequel, we will denote the first and second order terms of $A$ by $P_A:=A-\overline{W_{L}}$.

\begin{Remark} Note that since the vector fields $L^{\alpha\beta}$ generate orthogonal transformations of the fibers and since the eigenfunctions $U_0$ corresponding to the lowest eigenvalue $\lambda_0 > 0$ of the Dirichlet problem for the flat unit ball are orthogonally invariant, we have
\begin{equation}\label{annulator}
P_{A}u = 0
\end{equation}
for all $u = u_0\,\overline{v}\in E_0$. This fact will be used frequently in the sequel. Compare \cite{pap1}, Lemma 7.
\end{Remark}

\subsection{A uniform regularity result for the induced metric}

In this and in the following subsection, we will prove a basic result about uniform regularity of the operator family
$\Delta (\eps)$, $\eps > 0$, which implies that the 2-Sobolev norms of the solutions $w(\eps)$ of $\Delta (\eps)w(\eps)
= f\in L^2(L(1),g_0)$ -- which by elliptic regularity exist and are contained in $\bsob^1\cap\sob^2(L(1),g_0)$ for each
individual $\eps > 0$ -- are in fact uniformly bounded for all $\eps > 0$.

To be precise, we will prove the following statement.

\begin{Proposition}\label{uniform1} Let $u\in \bsob^1\cap\sob^2 (L(1),g_0)$. Then there is a constant $K_1 > 0$ ans some $\eps_1 > 0$ such
that
\begin{equation*}
\Vert u \Vert_2 \leq K_1 \left( \Vert \Delta(\eps) u \Vert_0 + \Vert u \Vert_0 \right)
\end{equation*} uniformly in $\eps < \eps_1$.
\end{Proposition}

The proof of the proposition splits up into a series of lemmas where $\Vert \Delta(\eps)u\Vert_0$ is decomposed into
different parts which can be related by Kato-type inequalities to the corresponding norm of the reference family. First,
we consider once again the reference family as the dominant term in the perturbation expansion of the induced family.

\begin{Lemma}\label{unterab} Let $u\in \bsob^1\cap\sob^2 (L(1),g_0)$ and $u^{\perp} := (1-E_0)u$. Then
\begin{equation*}
\Vert\Delta_0(\eps)u\Vert_0^2 \geq \frac{1}{\eps^2}\Vert\Delta_0^Vu^{\perp}\Vert^2_0 + \Vert\Delta_0^H u\Vert^2_0 +
\frac{2}{\eps}\langle \Delta_0^Hu^{\perp},\Delta_0^Vu^{\perp} \rangle_0
\end{equation*} and we have $\langle \Delta_0^Hu^{\perp},\Delta_0^Vu^{\perp} \rangle_0\geq 0$.
\end{Lemma}

\begin{Proof} Since $u$ is contained in the domain of the Dirichlet operators, we may use its spectral decomposition and
obtain
\begin{eqnarray*}
&& \Vert\Delta_0(\eps)u\Vert^2_0 \\
&=& \frac{1}{\eps^4}\Vert(\Delta_0^V-\lambda_0)u\Vert_0^2 + \Vert \Delta_0^H u\Vert^2_0 + \frac{2}{\eps^2}\langle
\Delta_0^Hu,(\Delta_0^V-\lambda_0)u \rangle_0 \\
&=& \frac{1}{\eps^2}\sum_{k > 0}\left\lbrack\frac{\lambda_k - \lambda_0}{\eps}\right\rbrack^2 \Vert E_{k}u\Vert_0^2 +
\Vert \Delta_0^H u\Vert^2_0 + \frac{2}{\eps}\sum_{k > 0}\frac{\lambda_k - \lambda_0}{\eps}\langle E_{k} u,\Delta_0^H u
\rangle_0.
\end{eqnarray*}
Now by
\begin{enumerate}
\item $\lbrack E_{k},\Delta_0^H\rbrack = 0$, $E_{k}^2 = E_{k}$, $E_{k}^+ = E_{k}$ which implies $$
\langle E_{k} u,\Delta_0^H u \rangle_0=\langle E_{k} u, \Delta_0^H E_{k} u \rangle_0\geq 0,$$
\item $\sum_{k > 0} E_{k} = 1 - E_0$, $u^{\perp} := (1-E_0)u$,
\item $\lambda_k - \lambda_0 / \eps > \lambda_k$ for $\eps < 1 - \lambda_0/\lambda_1$ (Lemma \ref{smoothEV}),
\end{enumerate} we obtain
\begin{eqnarray*}
\Vert\Delta_0(\eps)u\Vert^2_0 &\geq& \frac{1}{\eps^2}\sum_{k > 0}\lambda_k^2 \Vert E_{k}u\Vert_0^2 + \Vert \Delta_0^H
u\Vert^2_0 + \frac{2}{\eps}\sum_{k > 0}\lambda_k \langle E_{k} u,\Delta_0^H u \rangle_0\\
&=& \frac{1}{\eps^2}\Vert\Delta_0^V u^{\perp}\Vert_0^2 + \Vert \Delta_0^H u\Vert^2_0 + \frac{2}{\eps} \langle \Delta_0^V
u^{\perp},\Delta_0^H u^{\perp} \rangle_0 .
\end{eqnarray*}
By $u^{\perp} := (1-E_0)u\in\bsob^1\cap\sob^2 (L(1),g_0)$ we may conclude from {\rm (i)}, Lemma \ref{commutative} and (\ref{horipositiv}) the second statement
\begin{eqnarray*}
\langle \Delta_0^V u^{\perp}, \Delta_0^H u^{\perp}\rangle_0 &=& \sum_{k \geq 0} \lambda_k\langle E_{k} u^{\perp},
\Delta_0^H u^{\perp}\rangle_0 \\ &=& \sum_{k \geq 0} \lambda_k\langle E_{k} u^{\perp},
\Delta_0^H E_ku^{\perp}\rangle_0 = \tau(E_ku^{\perp})\geq 0 .
\end{eqnarray*}
\end{Proof}

That implies the following individual inequalities.

\begin{Lemma}\label{Abschatz} For $u\in \bsob^1\cap\sob^2 (L(1),g_0)$, we have $$
\Vert \Delta_0 (\eps) u \Vert_0 \geq\left\lbrace
\begin{array}{c} \frac{1}{\eps}\Vert \Delta_0^V u^{\perp}\Vert_0 \\
\Vert \Delta_0^H u\Vert_0 \\
\sqrt{\frac{2}{\eps}} \langle \Delta_0^V u^{\perp},\Delta_0^H u^{\perp} \rangle_0^{1/2}
\end{array}\right. . $$
\end{Lemma}

Recall that $u^{\perp} := (1-E_0)u$ and note that we use the same local Fermi coordinates and
notations as in Theorem \ref{Fermi Metrik} and that we always sum
over indices occuring twice if not stated otherwise.

In order to estimate the norm $\Vert - \Vert_{2k}$ of functions
$u\in\sob^{2k}(L(1),g_0)$ in the 2k-Sobolev spaces, we will first
localize the situation by a {\em basic} partition of unity and
explain how the bounds in local coordinates
serve to establish a global bound.

Let thus $\UU := U_{\iota , \iota\in I}$ be a finite open covering
of the submanifold $L$ which is trivializing for the normal
bundle, in particular $\pi^{-1}(U_{\iota}) \equiv U_{\iota} \times
B$, where $B\subset\R^{m-l}$ is again the euclidean unit ball. Let
furthermore $\chi_{\iota,\iota\in I}$ be a subordinated smooth
partition of unity on $L$ where we assume for convenience that
$\chi_{\iota}\geq 0$ for all $\iota\in I$. Then, the functions
$\overline{\chi_{\iota}}:=\chi_{\iota}\circ\pi$ form a smooth
partition of unity on $L(1)$ consisting only of basic functions.
Local coordinates thus yield an embedding $$ \pi^{-1}(U_{\iota})
\equiv U_{\iota} \times B\equiv V_{\iota}\subset \R^{l}\times
\R^{m-l} $$ and we can thus think of $V_{\iota}\equiv
\pi^{-1}(U_{\iota})$ as a regular subset of $\R^m$ and of
$\overline{\chi_{\iota}}u$ as a real function.

The next result is another Kato-type estimate, in this case for the perturbing operator.

\begin{Lemma}\label{nochnkato} Let $u\in\bsob^1\cap\sob^2(L(1),g_0)$. Then there is some uniform $C>0$ such that $$
\Vert Au\Vert_0 \leq C\,\left( \eps\Vert \Delta_0(\eps) u\Vert_0 + \Vert u\Vert_0\right) $$
for all $\eps > 0$.
\end{Lemma}

\begin{Proof} By (\ref{Alokal}), by (\ref{annulator}), the fact that the quadratic form $\langle
\d u, \d u\rangle_{H} = \langle u, P_{A}u\rangle_{0}$ vanishes for $u\in E_{0}$, and by the boundedness of
$\overline{W_{L}}$, we conclude that there is some $C_{W}> 0$ such that
\begin{eqnarray*}
\Vert Au\Vert_{0} &=&\Vert (P_{A} + \overline{W_{L}})u \Vert_{0} =  \Vert P_{A}u^{\perp} + \overline{W_{L}}\,u
\Vert_{0}\\
&\leq& \Vert P_{A}u^{\perp}\Vert_{0} + C_{W}\Vert u \Vert_{0}.
\end{eqnarray*}
Now we use the basic partition of unity $\overline{\chi}_{\iota,\iota\in I}$. Since the $\overline{\chi}_{\iota}$
are constant along the fibers, we have $E_{0}(\overline{\chi}_{\iota}u)=\overline{\chi}_{\iota}E_{0}u$ and hence $$
(\overline{\chi}_{\iota}u)^{\perp} = \overline{\chi}_{\iota}u -
E_{0}(\overline{\chi}_{\iota}u)=\overline{\chi}_{\iota}(u-E_{0}u)=\overline{\chi}_{\iota}u^{\perp}. $$
Thus, using the local expression for $P_A$ and the smoothness and boundedness of the coefficients, there are
numbers $C_{\iota,\iota\in I} > 0$ such that
\begin{eqnarray*}
\Vert P_{A}u^{\perp}\Vert_{0}  &\leq& \sum_{\iota\in I}\Vert
P_{A}(\overline{\chi}_{\iota}u)^{\perp}\Vert_{0} \\
&=& \sum_{\iota\in I}\Vert
\frac{1}{12}R_{\beta\nu\alpha\mu}L_{\beta\nu}L^{\alpha\mu} (\overline{\chi}_{\iota}u)^{\perp}\Vert_{0} \\
&\leq& \sum_{\iota\in I}C_{\iota} \Vert \overline{\chi}_{\iota} L_{\beta\nu}L_{\alpha\mu} u^{\perp}\Vert_{0} \\
\end{eqnarray*}
The vector fields $L^{\alpha\beta}$ are well-known to generate orthogonal transformations of the fibers. Hence the
operators $L^{\beta\nu}L^{\alpha\mu}$ are second order differential operators with with the same direct integral
structure as the vertical Laplacian $\Delta_0^V$. Since $(\overline{\chi}_{\iota}u)^{\perp}\in
\bsob^1\cap\sob^2(\pi^{-1}(q))$ and by elliptic regularity of the Dirichlet Laplacian on the flat unit ball, we
may use the vertical Laplacian to construct a Sobolev-2-norm on the individual fibers and we obtain that there are
numbers $C^{\beta\nu\alpha\mu}_{\iota,1},C^{\beta\nu\alpha\mu}_{\iota,2},C^{\beta\nu\alpha\mu}_{\iota,3} > 0$ such that
\begin{eqnarray*}
&& \Vert \overline{\chi}_{\iota}L_{\beta\nu}L_{\alpha\mu} u^{\perp}\Vert_{0}^2=\int_{L}\m_{L}(dq)\,\overline{\chi}_{\iota}^2\Vert L_{\beta\nu}L_{\alpha\mu} u^{\perp}_{q}\Vert^2_q \\
&\leq& C^{\beta\nu\alpha\mu}_{\iota,1}\int_{L}\m_{L}(dq)\,\Vert
u^{\perp}_{q}\Vert_{\sob^2(\pi^{-1}(q))}^2 \\
&\leq& C^{\beta\nu\alpha\mu}_{\iota,2}\int_{L}\m_{L}(dq)\,\left(\Vert
u^{\perp}_{q}\Vert_q + \Vert
D_q u^{\perp}_{q}\Vert_q\right)^2 \\
&=& C^{\beta\nu\alpha\mu}_{\iota,2}\left(\Vert u^{\perp}\Vert_{0} + \Vert
\Delta_0^V u^{\perp}\Vert_{0}\right)^2 .
\end{eqnarray*}
Hence by Lemma \ref{Abschatz}
\begin{equation*}
\label{een}
\Vert L^{\beta\nu}L^{\alpha\mu} (\overline{\chi}_{\iota}u)^{\perp}\Vert_{0}
\leq C^{\beta\nu\alpha\mu}_{\iota,2}\left(\Vert u^{\perp}\Vert_{0} + \eps\,\Vert
\Delta_{0}(\eps)u\Vert_{0}\right) .
\end{equation*}
Inserting this into the estimate for $\Vert P_Au^{\perp}\Vert_0$ above together with the
treatment of the zero order term yields the statement.
\end{Proof}

\subsection{Proof of Proposition \ref{uniform1}}

\begin{Proof} As explained above $$
\Delta (\eps) = \Delta_{0}(\eps) + A + \eps R(\eps) $$
where $A$ and $R(\eps)$ are second order differential operators with smooth coefficients. For $A$ and $R(\eps)$ these
coefficients together with all their derivatives are uniformly bounded. With some suitable $b > 0$ we obtain by Lemma \ref{nochnkato} and since $B(\eps)$ are second order differential operators with uniformly bounded coefficients
\begin{eqnarray*}
\Vert\Delta(\eps)u\Vert_0 + b\,\Vert u\Vert_0 &=& \Vert \Delta_{0}(\eps)u + Au + \eps R(\eps)u\Vert_0 + b\,\Vert u\Vert_0 \\
&\geq& \Vert \Delta_{0}(\eps)u\Vert_0 -  \Vert Au + \eps R(\eps)u\Vert_0 + b\,\Vert u\Vert_0 \\
&\geq& \Vert \Delta_{0}(\eps)u\Vert_0 -   \Vert Au\Vert_0 - \eps \Vert R(\eps)u\Vert_0 + b\,\Vert u\Vert_0\\
&\geq& \Vert \Delta_{0}(\eps)u\Vert_0 -  C( \eps\Vert \Delta_0(\eps)u\Vert_0 + \Vert u\Vert_0) - \eps d\Vert u\Vert_2 + b\,\Vert u\Vert_0 \\
&=& (1 - C\eps)\,\Vert \Delta_{0}(\eps)u\Vert_0 + (b- C)\, \Vert u\Vert_0 - \eps d\Vert u\Vert_2 .
\end{eqnarray*}
By Lemma \ref{normlap}, $\Vert u\Vert_2 \leq d' (\Vert\Delta_0u\Vert_0 + \Vert u\Vert_0)$ and thus
$$
\Vert\Delta(\eps)u\Vert_0 + b\,\Vert u\Vert_0
\geq (1 - C\eps )\,\Vert \Delta_{0}(\eps)u\Vert_0 + (b- C-dd'\eps)\, \Vert u\Vert_0) - dd'\eps\,\Vert\Delta_0 u\Vert .
$$
By Lemma \ref{unterab}, we have (recall that $\eps \leq 1$)
\begin{eqnarray*}
\Vert\Delta(\eps)u\Vert_0 &\geq& \Vert \Delta_0^Vu^{\perp} + \Delta_0^H u\Vert_0 \\
&=& \Vert \Delta_0^Vu^{\perp} + \lambda_0 E_0u + \Delta_0^H u - \lambda_0E_0u\Vert_0 \\
&\geq& \Vert \Delta_0u \Vert - \lambda_0\Vert u\Vert_0 .
\end{eqnarray*}
Now we choose $\eps > 0$ sufficiently small and $b >0$ sufficiently large such that $c_1 = 1 - (C+dd')\eps, c_2 = b- C-dd'\eps - (1-C\eps)\lambda_0 > 0$. That implies
$$
\Vert\Delta(\eps)u\Vert_0 + b\,\Vert u\Vert_0
\geq c_1\,\Vert \Delta_{0}u\Vert_0 + c_2\, \Vert u\Vert_0
$$
and by Lemma \ref{normlap} again, the expression on the right hand side is equivalent to the 2-Sobolev norm on $\bsob^1\cap\sob^2(L(1),g_0)$. That implies the statement.
\end{Proof}

\subsection{Boundary conditions for ${\mathrm{C^{\infty}}}$-vectors of $\Delta (\eps)$}

Still, our primary aim is to uniformly estimate the Sobolev norms of the family $u(\eps) := e^{ -\frac{t}{2}\Delta
(\eps)}u_{\eps}$ for some fixed value $t > 0$. We will prove that for $C^{\infty}$-vectors, the $2k$-Sobolev norms can be defined using $k$th powers of the Laplacian $\Delta_{0}$ associated to the reference metric. For $k=1$, we have
$\DD(\Delta (\eps)) = \DD(\Delta_{0})=\bsob^1\cap\sob^2(L(1),g_{0})$ where $\Delta_{0}$ here denotes the corresponding
Dirichlet operator. This justifies the use of the spectral theorem for $\Delta (\eps)$ and $\Delta_{0}$ in the preceding subsections. However, for $k>1$, the domains do no longer coincide and we have to be more careful about boundary
conditions. Since $u(\eps)$ is a $C^{\infty}$-vector for $\Delta (\eps)$ ({\em even analytic}, see i.e.
\cite{ReeSim:75b}, X.6), we will concentrate on boundary conditions for $C^{\infty}$-vectors and assume therefore
throughout this subsection that $u\in C^{\infty}(\Delta(\eps)):=\bigcap_{k\geq 1}\DD\left(
\Delta(\eps)^k\right)$. Note that by $\DD\left(
\Delta(\eps)^k\right)\subset \bsob^1\cap\sob^{2k}(L(1),g_{0})$ and by the {\em Sobolev embedding theorem}, the smooth
vectors can be considered as smooth functions and we may therefore simply skip a discussion of trace-maps in determining the boundary values.

In principle, we use the same idea as in the proof of Lemma \ref{normlap}, but now we have to seriously discuss the boundary conditions. We start with two lemmas which are verified in local coordinates.

\begin{Lemma}\label{Rand} Let $u\in C^{\infty}(\Delta(\eps))$. $u\bd=0$ implies $\Delta_0^Hu\bd = Au\bd=0$.
\end{Lemma}

\begin{Proof} In local coordinates (\ref{Alokal}), we see that $A$ contains only derivatives that generate rotations of the fiber and therefore leave the boundary conditions unaffected.
The zero order term of $A$ is bounded and continuous and thus we
have also $Au\bd=0$ if $u\bd =0$. For $\Delta_0^H$, we have
$$
\Delta_0^H = - \frac{1}{\sqrt{\det g_L}} X_i (g_L^{ij}\,\sqrt{\det g_L})X_j - g_L^{ij}X_i X_j
$$
where the $X_i$ are the horizontal vector fields (\ref{horizont}) which consist of differentiation in direction of the submanifold and of rotations of the fiber as in the case of $A$.
These derivatives do not affect zero boundary values at the boundary $\vert w\vert =1$ either.
\end{Proof}

\begin{Lemma}\label{Komm} Let $u\in C^{\infty}(\Delta(\eps))$. We have $\lbrack\Delta_0^V , \Delta_0^H \rbrack u =
\lbrack\Delta_0^V , A \rbrack u = 0$.
\end{Lemma}

\begin{Proof} In local coordinates, $\Delta_0^V$ consists only of derivatives in $w$-direction and will thus not affect
the coefficients of $A$ which are basic functions and thus constant along the fibers. Furthermore, it commutes with every generator of a rotation of the
fiber and hence with all $L^{\alpha\mu}$ since these vector fields generate isometries of the fiber. The statement follows thus from $\lbrack \Delta_0^V , A \rbrack u = \sum_{\iota\in I}\lbrack
\Delta_0^V , A \rbrack (\chi_{\iota}u) = 0$. $\lbrack\Delta_0^V , \Delta_0^H \rbrack u = 0$
was already discussed in Corollary \ref{commutative}.
\end{Proof}

The following proposition will be the main tool in the proof of the uniform estimate in the subsequent subsection. It
states that the boundary conditions defining $\DD(\Delta (\eps)^n)$ can in fact be represented as inhomogeneous boundary conditions for the regular elliptic boundary problem from Section \ref{Lapnorm}.

\begin{Proposition}\label{BigRand} Let $u\in C^{\infty}(\Delta(\eps))$. For $n\geq 1$ we have
\begin{equation*}
\Delta_0^nu\bd = \eps^3\, T_n(\eps)u\bd
\end{equation*}
where $T_n(\eps)$ is a family of differential operators of order $2n$ on $L(1)$ with smooth coefficients which are
bounded together with all their derivatives uniformly in $\eps > 0$.
\end{Proposition}

\begin{Proof} We prove this statement by induction. {\rm (i)} Let first be $n=1$. By the perturbation expansion for
$\Delta (\eps)$ we have $$
\Delta(\eps)u\bd = \left\lbrace \frac{\Delta_0^V - \lambda_0}{\eps^2}u + (\Delta_0^H + A)u + \eps \,
R(\eps)u\right\rbrace\bd $$
where $R(\eps)$ is a second order operator with uniformly bounded coefficients. By Lemma \ref{Rand}, that implies $$
\Delta(\eps)u\bd = \left\lbrace \frac{1}{\eps^2}\Delta_0^Vu  + \eps \, R(\eps)u\right\rbrace\bd $$
and since by assumption $\Delta(\eps)u\bd = 0$, we have $\Delta_0^Vu\bd  = - \eps^3 \, R(\eps)u\bd $ and thus again by
Lemma \ref{Rand} $$
\Delta_0 u\bd = \Delta_0^Vu\bd  + \Delta_0^Hu\bd = - \eps^3 R(\eps)u\bd. $$
Letting $T_1(\eps) := -R(\eps)$ yields thus the statement in the case $n=1$. {\rm (ii)} For the step from $n$ to $n+1$,
we have
\begin{eqnarray*}
0 &=& \Delta(\eps)^{n+1}u\bd = (\eps^2 \Delta(\eps))^{n+1}u\bd \\
&=& (\Delta_0^V + \eps^3 R(\eps)) (\eps^2
\Delta(\eps))^{n}u\bd \\
&=& \big\lbrace \Delta_0^V(\Delta_0^V - \lambda_0)^n + \eps^2 \sum_{s=0}^{n-1} \Delta_0^V(\Delta_0^V -
\lambda_0)^s(\Delta_0^H + A)(\Delta_0^V - \lambda_0)^{n-s-1} \\
&& + \eps^3 R_n(\eps)\big\rbrace u\bd \\
\end{eqnarray*}
where $R_n(\eps)$ is of order $2n+2$ with uniformly bounded smooth coefficients. By Lemma \ref{Komm}, we have thus
\begin{eqnarray*}
0 &=& \big\lbrace \Delta_0^V(\Delta_0^V - \lambda_0)^n + n\eps^2 (\Delta_0^H + A)\Delta_0^V(\Delta_0^V -
\lambda_0)^{n-1}  +
\eps^3 R_n(\eps)\big\rbrace u\bd \\
&=& \Delta_0^V(\Delta_0^V - \lambda_0)^n u\bd   + \eps^3 R_n(\eps) u\bd \\
&& + n\eps^2 \sum_{s=0}^{n-1}\left(
\begin{array}{c} n-1 \\
s
\end{array}\right)(-\lambda_0)^{n-s-1}(\Delta_0^H + A)(\Delta_0^V)^{s+1}u\bd \\
&=& \Delta_0^V(\Delta_0^V - \lambda_0)^n u\bd   + \eps^3 R_n(\eps) u\bd \\
&& + n\eps^2 \sum_{s=0}^{n-1}\left(
\begin{array}{c} n-1 \\
s
\end{array}\right)(-\lambda_0)^{n-s-1}(\Delta_0^H + A)\left\lbrace(\Delta_0^V)^{s+1}u - \eps^3 T_{s+1}(\eps)\right\rbrace
u\bd \\
&& + n\eps^2 \sum_{s=0}^{n-1}\left(
\begin{array}{c} n-1 \\
s
\end{array}\right)(-\lambda_0)^{n-s-1}(\Delta_0^H + A)\eps^3 T_{s+1}(\eps)u\bd .
\end{eqnarray*}
By induction hypothesis (note that $s+1\leq n$), we have $$
((\Delta_0^V)^{s+1}u - \eps^3 T_{s+1}(\eps)) u\bd =0 $$
and hence again by Lemma \ref{Rand}
\begin{eqnarray*}
0 &=& \Delta_0^V(\Delta_0^V - \lambda_0)^n u\bd   + \eps^3 R_n(\eps) u\bd \\
&& + \eps^5 n \sum_{s=0}^{n-1}\left(
\begin{array}{c} n-1 \\
s
\end{array}\right)(-\lambda_0)^{n-s-1}(\Delta_0^H + A) T_{s+1}(\eps)u\bd \\
&=& \Delta_0^V(\Delta_0^V - \lambda_0)^n u\bd   + \eps^3 R_n(\eps) u\bd \\
&& + \eps^5 n \sum_{s=0}^{n-1}\left(
\begin{array}{c} n-1 \\
s
\end{array}\right)(-\lambda_0)^{n-s-1}(\Delta_0^H + A) T_{s+1}(\eps)u\bd \\
&=& \Delta_0^V(\Delta_0^V - \lambda_0)^n u\bd   + \eps^3 (R_n(\eps) + \eps^2\tilde{R}_n(\eps))u\bd .
\end{eqnarray*}
Thus, since $R_n(\eps) + \eps^2\tilde{R}_n(\eps)$ is a family of differential operators of order $2n + 2$ with smooth
and uniformly bounded coefficients as required, we obtain $$
\Delta_0^V(\Delta_0^V - \lambda_0)^n u\bd = \eps^3 \tilde{T}_{n+1}(\eps)u\bd . $$
That implies $$
(\Delta_0^V)^{n+1} u\bd = \eps^3 \tilde{T}_{n+1}(\eps)u\bd -\sum_{s=0}^{n-1}\left( \begin{array}{c} n \\
s
\end{array} \right)(-\lambda_0)^{n-s}(\Delta_0^V)^{s+1}u\bd $$
and another application of the induction hypothesis to the sum on the right hand side yields (again $s+1\leq n$) the
statement $$
(\Delta_0^V)^{n+1} u\bd = \eps^3 T_{n+1}(\eps)u\bd . $$
\end{Proof}

This characterization of the boundary values enables us to proof the analogue of Lemma \ref{normlap}
for the subset
$$
\DD(\eps_0,k) := \cup_{\eps\leq\eps_0}\DD (\Delta (\eps)^k)\subset \sob^{2k}(L(1),g_0)
$$
which contains the families of semigroups generated by the rescaled Laplacians associated to the induced metric. For $\Delta_0(\eps)$, the domains $\DD (\Delta_0 (\eps)^k)$ are independent of $\eps$. The complication is that we now have to deal with the fact that for the domains of powers of $\Delta (\eps)$ this is no longer true in general. But we can control the differences of the domains by representing the functions as solutions of the elliptic boundary problem from Section \ref{Lapnorm} where now the boundary conditions {\em depend on $\eps$}. By the preceding proposition, we will be able to control these boundary conditions in a suitable way.

\begin{Proposition}\label{LSindu} There is some $\eps_0>0$ such that the $2k$-Sobolev norm is equivalent to the norm $$
\vert\Vert u \Vert\vert_{2k} := \Vert u \Vert_0 + \Vert \Delta_0^k u \Vert_0 $$
on $\DD(\eps_0,k) = \cup_{\eps\leq\eps_0}\DD (\Delta (\eps)^k)\subset \sob^{2k}(L(1),g_0)$.
\end{Proposition}

\begin{Proof} Let $u\in\DD(\Delta(\eps)^k)$. It is again sufficient to prove that for some $C>0$ that does not depend on
$\eps$ we have the inequality $$
\Vert u\Vert_{2k} \leq C\left(\Vert u \Vert_0 + \Vert \Delta_0^k u \Vert_0\right). $$
First of all, we use the fact that a function $u\in\DD(\Delta(\eps)^k)$ solves the boundary value problem from Section
\ref{Lapnorm}, but with different boundary conditions. In particular, set $f(u):=\Delta_0^k u\in L^2(L(1),g_0)$ and
\begin{eqnarray*}
g_1 &:=& u\bd = 0, \\
g_2&:=&\Delta_0u\bd=\eps^3\,T_2(\eps)u\bd , \\
&...& \\
g_{k} &:=&\Delta_0^{k-1}u\bd=\eps^3\,T_k(\eps)u\bd .
\end{eqnarray*}
Using the notation from Section \ref{Lapnorm} again, we see that $u$ solves the regular elliptic boundary problem $$
\LL u = f(u),\BB_ju=g_j, $$
$j=1,...,k$, with inhomogeneous boundary conditions. But that implies (see e.g. \cite{Agr:97}, Theorem 2.2.1, p. 16)
that
\begin{equation}
\label{hierundjetzt}
\Vert u\Vert_{2k} \leq C\left(\Vert u \Vert_0 + \Vert \Delta_0^k u \Vert_0 + \sum_{j=2}^k\Vert
g_j\Vert_{\sob^{2(k-j)-1/2}(
\partial L(1),g_0)}\right).
\end{equation}
It is important to note that the constant $C>0$ only depends on the boundary problem and therefore not on $\eps > 0$.
The only part here that still depends crucially on $\eps$ and that has to be taken care of are the {\em boundary values}
$g_j$. To do so, we note that the operators $T_j(\eps)$ are of order $2j$ with smooth, uniformly bounded coefficients
and they are defined on the whole tube. Hence, by continuity of the trace map
$\sob^{s}(L(1),g_0)\to\sob^{s-1/2}(\partial L(1),g_0)$, we have constants $c_2,...,c_k > 0$ not depending on $\eps$ such
that $$
\Vert g_j\Vert_{\sob^{2(k-j)-1/2}(
\partial L(1),g_0)}=\eps^3\Vert T_j(\eps)u\bd\Vert_{\sob^{2(k-j)-1/2}(\partial L(1),g_0)} \leq \eps^3 c_j\,\Vert
u\Vert_{2k}. $$
For all $\eps$ less than some $\eps_0 > 0$, we can thus absorb these terms into the left hand side of the inequality
(\ref{hierundjetzt}) above. That implies the statement.
\end{Proof}

\subsection{The proof of Proposition \ref{compactness}}\label{comconv}

Now we come to the last step in the proof of Theorem \ref{MainSemi2}, the proof of Proposition \ref{compactness}.

First of all, we state an interpolation result that we will have to use in the sequel. It is a global version of the
local statement from \cite{Agm:65}, Theorem 3.4, p. 25.

\begin{Lemma}\label{interpol}{\bf (Interpolation)} Let $1\leq k \leq 2n+1$ be integer. Then, for all $1 \geq \alpha > 0$
there are constants $C(k,\alpha)>0$ and $C(k)> 0$ such that $$
\Vert u \Vert_{k} \leq C(k)\,\left(\alpha \Vert u\Vert_{2n+2} + C(k,\alpha)\Vert u\Vert_0\right). $$
\end{Lemma}

\begin{Proof} By \cite{Tay:99}, 4.4, p. 284 ff., the spaces $\sob^k(L(1),g_{0})$ with $k$ as above are
interpolation spaces of the form $$
\sob^k (L(1),g_{0}) = \lbrack L^2(L(1),g_{0}),\sob^{2n}(L(1),g_{0})\rbrack_{k/2n}. $$
By the {\em Calderon-Lions interpolation theorem} (see \cite{ReeSim:75b}, Theorem IX.20, p. 37) that implies the
estimate $$
\Vert u \Vert_{k}\leq \Vert u\Vert_{0}^{1-k/n} \,\Vert u \Vert_{2n}^{k/n}. $$
Using {\em Young's inequality} $ab\leq a^p/p + b^q/q$, for $a,b\geq 0$, $p,q\geq 1$, $1/p + 1/q =1$, we thus obtain $$
\Vert u \Vert_{k}\leq \frac{1}{\alpha}\Vert u\Vert_{0}^{1-k/2n} \,\alpha\Vert u \Vert_{2n}^{k/2n}\leq
\frac{1-k/2n}{\alpha^{1/(1-k/2n)}}\Vert u\Vert_{0} + \frac{2n}{k}\alpha^{2n/k}\Vert u \Vert_{2n} $$
where $p = 1/(1-k/2n)$, $q=2n/k$. By $2n/k > 1$, $\alpha^{2n/k}\leq \alpha$ and the statement is proved.
\end{Proof}

Now we prove a uniform regularity result for the induced family on the subspace of the associated boundary problem.

\begin{Proposition}\label{uniform2} Let $u\in\bigcap_{k\geq 1}\DD\left( \Delta(\eps)^k\right)$. Then, for all $n\geq 0$
there are numbers $K_n > 0$, $\eps_n > 0$ such that $$
\Vert u \Vert_{2n+2} \leq K_n \left( \Vert \Delta (\eps) u\Vert_{2n} + \Vert u\Vert_0\right) $$
uniformly for all $\eps \leq \eps_n$.
\end{Proposition}

\begin{Proof} We will prove this statement by induction. The case $n=1$ is given by Proposition \ref{uniform1}. The step
from $n$ to $n+1$ is provided by the following chain of arguments: We can apply the equivalent description of the
Sobolev norms by powers of $\Delta_{0}$ from Proposition \ref{LSindu} and obtain
\begin{eqnarray*}
\Vert u \Vert_{2n+2} &\leq& C \left( \Vert u \Vert_{0} + \Vert \Delta_0^{n+1} u\Vert_{0}\right) \\
&\leq & C \left( \Vert u \Vert_{0} + \Vert \Delta_0^{n} u\Vert_{2}\right) \\
&\leq & C \left( \Vert u \Vert_{0} + \Vert \Delta_0^n u - \eps^3 T_n(\eps)u\Vert_{2} + \eps^3 \Vert
T_n(\eps)u\Vert_2\right)
\end{eqnarray*}
where we insert the operators $T_n(\eps)$ from Proposition \ref{BigRand} which were globally defined on $L(1)$. Note that by Proposition \ref{BigRand}, we have $$
\Delta_0^n u - \eps^3 T_n(\eps)u \in\bsob^1\cap\sob^2 (L(1),g_0). $$
and we may apply Proposition \ref{uniform1} to obtain
\begin{eqnarray*}
\Vert u \Vert_{2n+2} &\leq & C  (\Vert u \Vert_{0}  + C\eps^3 \Vert T_n(\eps)u\Vert_2 ) \\
&& + CK_1( \Vert \Delta (\eps)(\Delta_0^n u - \eps^3 T_n(\eps)u)\Vert_0 + \Vert \Delta_0^n u - \eps^3
T_n(\eps)u\Vert_0)\\
&\leq & C  (\Vert u \Vert_{0}  + a\eps^3 \Vert u\Vert_{2n+2} ) \\
&& + CK_1( \Vert \Delta (\eps)(\Delta_0^n u - \eps^3 T_n(\eps)u)\Vert_0 + \Vert \Delta_0^n u - \eps^3
T_n(\eps)u\Vert_0)\\
&\leq & \tilde{C}_1  (\Vert u \Vert_{0}  + \eps^3 \Vert u\Vert_{2n+2} + \Vert u \Vert_{2n}) \\
&& + CK_1(\Vert \Delta (\eps)\Delta_0^n u \Vert_0 + \eps \Vert \eps^2\Delta(\eps) T_n(\eps)u\Vert_0) \\
&\leq & \tilde{C}_2  (\Vert u \Vert_{0}  + (\eps^3 + \eps) \Vert u\Vert_{2n+2} + \Vert u \Vert_{2n}) \\
&& + CK_1\,\Vert \Delta (\eps)\Delta_0^n u \Vert_0  \\
\end{eqnarray*}
where we use the fact that $\eps^2\Delta(\eps)$ is also a differential operator of second order with uniformly bounded
smooth coefficients. By Lemma \ref{interpol}, we see that for $\eps$ sufficiently small, we may absorb the $\Vert - \Vert_{2n}$-term and the $\eps\Vert -\Vert_{2n + 2}$-term into the $\Vert -
\Vert_{2n+2}$-term on the left hand side of the inequality and that we obtain after changing the constant if necessary
$$
\Vert u \Vert_{2n+2}\leq  \tilde{C}_3  (\Vert u \Vert_{0}   + \Vert \Delta (\eps)\Delta_0^n u \Vert_0 ). $$
But now by Lemma \ref{Komm} and Proposition \ref{DecoReference}, we have $\lbrack \Delta_0(\eps),\Delta_0\rbrack = 0$ an
we get $$
\lbrack \Delta(\eps),\Delta_0^n\rbrack = \lbrack \Delta_0(\eps),\Delta_0^n\rbrack + \lbrack A + \eps R(\eps),
\Delta_0^n\rbrack = \lbrack A + \eps R(\eps), \Delta_0^n\rbrack $$
and this is a differential operator with smooth and uniformly bounded coefficients of order $2n+1$. Hence
\begin{eqnarray*}
\Vert u \Vert_{2n+2} &\leq & \tilde{C}_3  (\Vert u \Vert_{0}   + \Vert \Delta (\eps)\Delta_0^n u \Vert_0 ) \\
&\leq & \tilde{C}_3  (\Vert u \Vert_{0} + \Vert \Delta_0^n\Delta (\eps) u \Vert_0  + \Vert \lbrack\Delta
(\eps),\Delta_0^n\rbrack u \Vert_0 ) \\
&\leq & \tilde{C}_4  (\Vert u \Vert_{0} + \Vert \Delta_0^n\Delta (\eps) u \Vert_0  + \Vert  u \Vert_{2n+1} ) \\
&\leq & \tilde{C}_4  (\Vert u \Vert_{0} + \Vert\Delta(\eps)u\Vert_0 + \Vert \Delta_0^n\Delta (\eps) u \Vert_0  + \Vert
u \Vert_{2n+1} ) .
\end{eqnarray*}
Absorbing again the $\Vert - \Vert_{2n+1}$-norm by Lemma \ref{interpol} in the same way as above finally yields $$
\Vert u \Vert_{2n+2}
\leq K_n (\Vert u \Vert_{0} + \Vert\Delta(\eps)u\Vert_{2n}) . $$
\end{Proof}

Let now $u_{\eps}$, $\eps > 0$, a sequence of initial conditions $u_{\eps}\in L^2(L(1),g_0)$ converging to $u\in L^2(L(1),g_0)$. Let $t > 0$ be fixed and again
\begin{equation*}
u(\eps, t ) := \exp \left( -\frac{t}{2}\Delta (\eps)\right)u_{\eps} \in \bigcap_{k\geq 1}\DD (\Delta(\eps)^k).
\end{equation*}

\begin{Lemma}\label{nochnenbound} For all $n\geq 1$ there are numbers $D_n>0$, $\eps_n > 0$ such that $$
\Vert u (\eps , t)\Vert_{2n} \leq D_n \left( \sum_{k=0}^n\Vert \Delta (\eps)^k u (\eps ,t)\Vert_0\right) $$
uniformly for all $\eps \leq\eps_n$.
\end{Lemma}

\begin{Proof} Again, we proceed by induction. The case $n=1$ is proven by Proposition \ref{uniform1}. Let now the
statement be true for $n$. Then by Proposition \ref{uniform2}
\begin{eqnarray*}
\Vert u (\eps,t)\Vert_{2n+2} &\leq& K_{n+1} (\Vert \Delta (\eps)u(\eps,t)\Vert_{2n} + \Vert u(\eps,t) \Vert_0 )\\
&\leq& K_{n+1} (D_n \left( \sum_{k=0}^n\Vert \Delta (\eps)^k \Delta(\eps)u(\eps,t)\Vert_0\right)+ \Vert u(\eps,t) \Vert_0 )\\
&\leq& D_{n+1}  \left( \sum_{k=0}^{n+1}\Vert \Delta (\eps)^k u(\eps,t)\Vert_0\right).
\end{eqnarray*}
\end{Proof}

That implies that we finally obtain compactness in the appropriate in Sobolev spaces.

\noindent{\bf Proof of Proposition \ref{compactness}:} By the spectral theorem $$
\Delta(\eps)^ku(\eps ,t)= \sum_{s\geq 0} \lambda_{s}(\eps)^k e^{-\frac{t}{2}\lambda_{s}(\eps)}\langle
u_{s}(\eps),u_{\eps}\rangle_{0}\,u_{s}(\eps) $$
since the Dirichlet-Laplacians are self-adjoint with discrete spectrum $$
0 < \lambda_{1}(\eps) \leq \lambda_{2}(\eps)
\leq ...$$
and associated normalized eigenfunctions $u_{s}(\eps)$ for all $\eps > 0$. Thus $$
\Vert \Delta(\eps)^k u(\eps,t)\Vert_{0}^2 =
\sum_{s\geq 0} \lambda_{s}(\eps)^{2k} e^{-t\lambda_{s}(\eps)}\vert\langle u_{s}(\eps),u_{\eps}\rangle_{0}\vert^2 $$
and by $x^{2k}e^{-tx}\leq (2k/t)^{2k}e^{-2k}=:c_{k}$ for $x\geq 0$, we obtain
\begin{equation}
\label{unifoermchen}
\Vert \Delta(\eps)^k u(\eps,t)\Vert_{0}^2 \leq c_{k}\Vert u_{\eps}\Vert_{0}^2 \leq C_k .
\end{equation}
since the sequence $u_{\eps}$ converges in $L^2(L(1),g_0)$. Combining this with Lemma \ref{nochnenbound} yields statement {\rm (i)}. For {\rm (ii)}, we consider $t(K) := \inf K$. If $a>0$, we have for all $t',t'' \in K$ (even for all $t',t'' > t(K)$) that
$$
\left\vert e^{-at'} - e^{-at''}\right\vert^2 \leq a^2e^{-2at(K)} \,\vert t' - t''\vert .
$$
That implies by {\rm (ii)}
\begin{eqnarray*}
\Vert \Delta(\eps)^k (u(\eps,t') &-& u(\eps,t''))\Vert_{0}^2 \\
&=& \sum_{s\geq 0} \lambda_{s}(\eps)^{2k} \left\vert e^{-\frac{t'\lambda_{s}(\eps)}{2}} - e^{-\frac{t''\lambda_{s}(\eps)}{2}} \right\vert^2 \,\vert\langle u_{s}(\eps),u_{\eps}\rangle_{0}\vert^2 \\
&\leq& \vert t' - t''\vert\,\sum_{s\geq 0} \lambda_{s}(\eps)^{2k+2} e^{-t(K)\lambda_{s}(\eps)} \,\vert\langle u_{s}(\eps),u_{\eps}\rangle_{0}\vert^2 \\
&\leq& \vert t' - t''\vert\,C_{k+1}
\end{eqnarray*}
and therefore uniform equicontinuity of the family $\AA$.

\hfill $\blacksquare$

\bibliographystyle{plain}

\end{document}